\documentclass[11pt,a4paper,twoside]{article}

\newif\ifpdf
\ifx\pdfoutput\undefined
\pdffalse 		
\else
\pdfoutput=1 		
\pdftrue
\fi

\ifpdf
\usepackage[pdftex]{graphicx}
\usepackage[colorlinks=true,%
            breaklinks=true,%
            citecolor=blue,%
	    urlcolor=black,%
	    bookmarksopen=true,%
            pdftex=true,%
            hyperindex=true,%
	    ]%
            {hyperref}
\else
\usepackage{graphicx}
\usepackage[nativepdf,%
            breaklinks=true,%
	    colorlinks=true,%
            citecolor=blue,%
	    urlcolor=black,%
	    bookmarksopen=true,%
            hyperindex=true,%
	    ]%
            {hyperref}
\fi

\title{\bf Geometric and Physical Interpretation \\
       of Fractional Integration \\
	and Fractional Differentiation}

\author{\textbf{Igor Podlubny}\\
	Department of Informatics and Control Engineering\\
	BERG Faculty, Technical University of Kosice\\
	B. Nemcovej 3, 04200 Kosice, Slovak Republic \\
	e-mail: \textit{\href{mailto:Igor.Podlubny@tuke.sk}{Igor.Podlubny@tuke.sk}}}

\begin{document}

\maketitle

\begin{abstract}
\noindent
A solution to the more than 300-years old problem
of geometric and physical interpretation of fractional
integration and differentiation
(i.e., integration and differentiation of an arbitrary
real order) is suggested for the Riemann-Liouville fractional integration and differentiation, the Caputo fractional differentiation, the Riesz potential,
and the Feller potential.
It is also generalized for giving a new geometric and physical
interpretation of more general convolution integrals
of the Volterra type.
	Besides this, a new physical interpretation
is suggested for the Stieltjes integral.



\end{abstract}

\section{Introduction}
\markboth{\footnotesize \textsf{\qquad Report TUKE--10--2001 \hfill I. Podlubny}}%
         {\footnotesize \textsf{Geometric and Physical Interpretation
         	 of Fractional Integration and Differentiation}}

It is generally known that integer-order derivatives and integrals
have clear physical and geometric interpretations, which
significantly simplify their use for solving applied problems
in various fields of science.

However, in case of fractional-order integration
and differentiation, which represent a rapidly
growing field both in theory and in applications
to real-world problems, it is not so.
Since the appearance of the idea of differentiation and
integration of arbitrary (not necessary integer) order
there was not any acceptable geometric and physical
interpretation of these operations for more than 300 years.
The lack of these interpretations has been acknowledged
at the first international conference on the fractional calculus
in New Haven (USA) in 1974 by including it in the list of
open problems \cite{Ross-conf}.
The question was unanswered, and therefore repeated
at the subsequent conferences at the University of Strathclyde (UK)
in 1984 \cite{McBride-Roach}
and at the Nihon University (Tokyo, Japan) in 1989 \cite{Nishimoto-conf}.
The round-table discussion
\cite{TMSF-96-Kiryakova,TMSF-96-Gorenflo,TMSF-96-Mainardi}
at the conference on transform methods and special functions
in Varna (1996) showed that the problem was still
unsolved, and since that time the situation,
in fact, still did not change.

Fractional integration and fractional differentiation
are generalisations of notions of integer-order
integration and differentiation, and include $n$-th
derivatives and $n$-folded integrals
($n$ denotes an integer number) as particular cases.
Because of this, it would be ideal to have such
physical and geometric interpretations of
fractional-order operators, which will provide also
a link to known classical interpretations of
integer-order differentiation and integration.

Since the need for the aforementioned geometric and
physical interpretations is generally recognised,
several authors attempted to provide them.
Probably due mostly to linguistical reasons,
much effort have been devoted to trying to relate
\textit{fractional} integrals and derivatives,
on one side, and \textit{fractal geometry}, on the other
\cite[and others]{Nigmatullin,cookie-set,self-similar-set,Monsrefi-Torbati-Hammond}.
However, it has been clearly shown by R.~Rutman
\cite{Rutman-critical,Rutman-constructive}
that this approach is inconsistent.

Besides those ``fractal-oriented'' attempts,
some considerations regarding interpretation of
fractional integration and fractional differentiation
were presented in \cite{Monsrefi-Torbati-Hammond}.
However, those considerations are, in fact, only
a small collection of selected examples of applications
of fractional calculus, in which hereditary effects
and self-similarity are typical for the objects
modelled with the help of fractional calculus.
Although each particular
problem, to which fractional derivatives or/and fractional
integrals have been applied, can be considered
as a certain illustration of their meaning,
the paper \cite{Monsrefi-Torbati-Hammond} cannot
be considered as a definite answer to the posed question.

A different approach to geometric interpretation
of fractional integration and fractional differentiation,
based on the idea of the contact of $\alpha$-th order,
has been suggested by F.~Ben Adda \cite{Benadda-JFC,Benadda}.
However, it is difficult to speak about an acceptable
\emph{geometric interpretation}
if one cannot see any picture there.

Obviously, there is still a lack of geometric and physical
interpretation of fractional integration and
differentiation, which is comparable with the
simple interpretations of their integer-order
counterparts.

In this paper we present a new approach to solution
of this challenging old problem.

We start with introducing a simple and really geometric
interpretation of several types of fractional-order integration:
the left-sided and the right-sided Riemann--Liouville fractional
integration, the Riesz potential, and the Feller potential.

Based on this, a physical interpretation of the Riemann--Liouville
fractional integration is proposed in terms of
inhomogeneous and changing (non-static, dynamic) time
scale.
Moreover, on this way we give
a new physical interpretation of the Stieltjes integral.
We also try to persuade the readers that the suggested
physical interpretation of fractional integration
is in line with the current views on space--time in physics.
	We also suggest physical interpretation
for the Riemann-Liouville fractional differentiation
and for the Caputo fractional differentiation.
	Finally, we show that the suggested approach to geometric
interpretation of fractional integration can be used
for providing a new geometric and physical interpretation
for convolution integrals of the Volterra type.

\section{Geometric interpretation of fractional \\ integration:
         Shadows on the walls}

In this section we first give a geometric interpretation of
left-sided and right-sided Riemann--Liouville fractional integrals,
and then consider the Riesz potential.

\subsection{Left-sided Riemann--Liouville fractional integral}

Let us consider the left-sided Riemann--Liouville fractional integral
\cite{Podlubny-FDE-book,SKM} of order $\alpha$,
\begin{equation}\label{FI}
	_{0}I_{t}^{\alpha} f(t)
	 =  	\frac{1}{\Gamma (\alpha)}
		\int\limits_{0}^{t}
		f(\tau) (t-\tau)^{\alpha - 1}
		d\tau,
\end{equation}
and write it in the form
\begin{equation}\label{FI-as-SI}
	_{0}I_{t}^{\alpha} f(t) =
		\int\limits_{0}^{t}
		f(\tau)
		dg_{t}(\tau),
\end{equation}
\begin{equation}\label{eq:g-function}
	g_{t}(\tau) = 	\frac{1}{\Gamma (\alpha + 1)}
			\Bigl\{
			t^{\alpha} - (t-\tau)^{\alpha}
			\Bigr\}.
\end{equation}

The function $g_{t}(\tau)$ has an interesting scaling property.
Indeed, if we take $t_{1}=kt$ and $\tau_{1}=k\tau$, then
\begin{equation}
	g_{t_{1}}(\tau_{1})= g_{kt}(k\tau)=
	k^{\alpha} g_{t}(\tau).
\end{equation}

Now let us consider the integral (\ref{FI-as-SI}) for a fixed $t$.
Then it becomes simply a Stieltjes integral, and we can utilize
G.~L.~Bullock's idea \cite{Bullock}.

Let us take the axes $\tau$, $g$, and $f$.
In the plane $(\tau,g)$ we plot the function $g_{t}(\tau)$
for $0\leq\tau\leq t$. Along the obtained curve we ``build a fence''
of the varying height $f(\tau)$, so the top edge
of the ``fence'' is a three-dimensional line
$(\tau, g_{t}(\tau), f(\tau))$, $0\leq\tau\leq t$.

\pagebreak

\begin{figure}[p]
\centering
\includegraphics[width=0.7\textwidth]{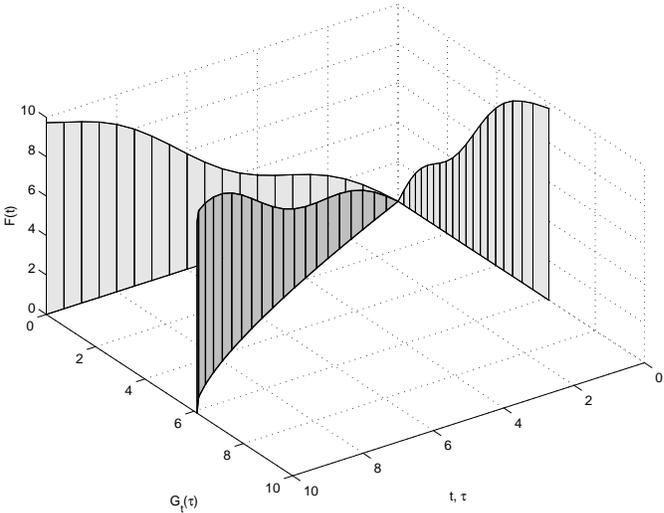}
\caption{\protect\parbox[t]{9cm}{The ``fence'' and its shadows:
$_{0}I_{t}^{1}f(t)$ and $_{0}I_{t}^{\alpha}f(t)$, \protect\\
for $\alpha =0.75$, $f(t)=t+0.5 \sin(t)$}}
\label{fig:fence}
\end{figure}

This ``fence'' can be projected onto two surfaces (see Fig.~\ref{fig:fence}):
\begin{itemize}
\item
the area of the projection of this
``fence'' onto the plane $(\tau, f)$ corresponds to
the value of the integral
\begin{equation}\label{one-fold-integral}
	_{0}I_{t}^{1}(t) = 	\int\limits_{0}^{t}
		f(\tau) d\tau;
\end{equation}
\item
the area of the projection of the same ``fence'' onto the plane $(g,f)$
corresponds to the value of the integral (\ref{FI-as-SI}), or,
what is the same, to the value of the fractional integral (\ref{FI}).
\end{itemize}

\begin{figure}[p]
\centering
\includegraphics[width=0.7\textwidth]{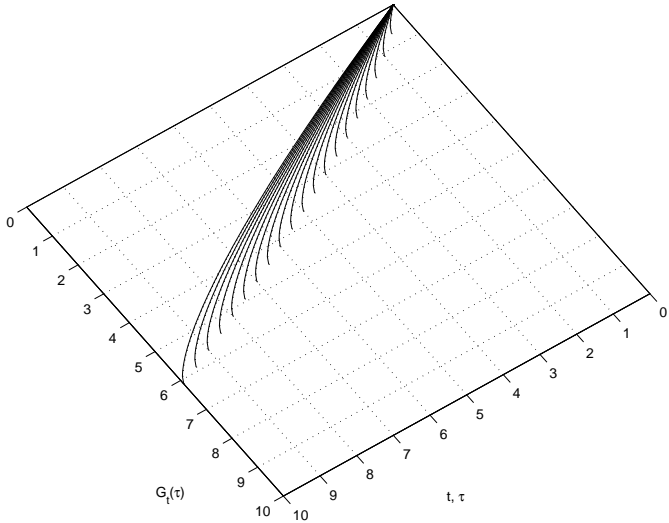}
\caption{\protect\parbox[t]{9cm}{The process of change of the fence
basis shape \protect\\
for $_{0}I_{t}^{\alpha}f(t)$, $\alpha=0.75$.}}
\label{fig:fencebas}
\end{figure}

In other words, our ``fence'' throws two shadows on two walls.
The first of them, that on the wall $(\tau,f)$, is the well-known
``area under the curve $f(\tau)$'', which is a standard geometric
interpretation of the integral (\ref{one-fold-integral}).
The ``shadow'' on the wall $(g,f)$ is a geometric interpretation
of the fractional integral (\ref{FI}) for a fixed $t$.

Obviously, for $g_{t}(\tau)=\tau$ both ``shadows'' are equal.
This shows that classical definite integration is a particular
case of the left-sided Riemann--Liouville
fractional integration even from the geometric point of view.

\begin{figure}
\centering
\includegraphics[width=0.7\textwidth]{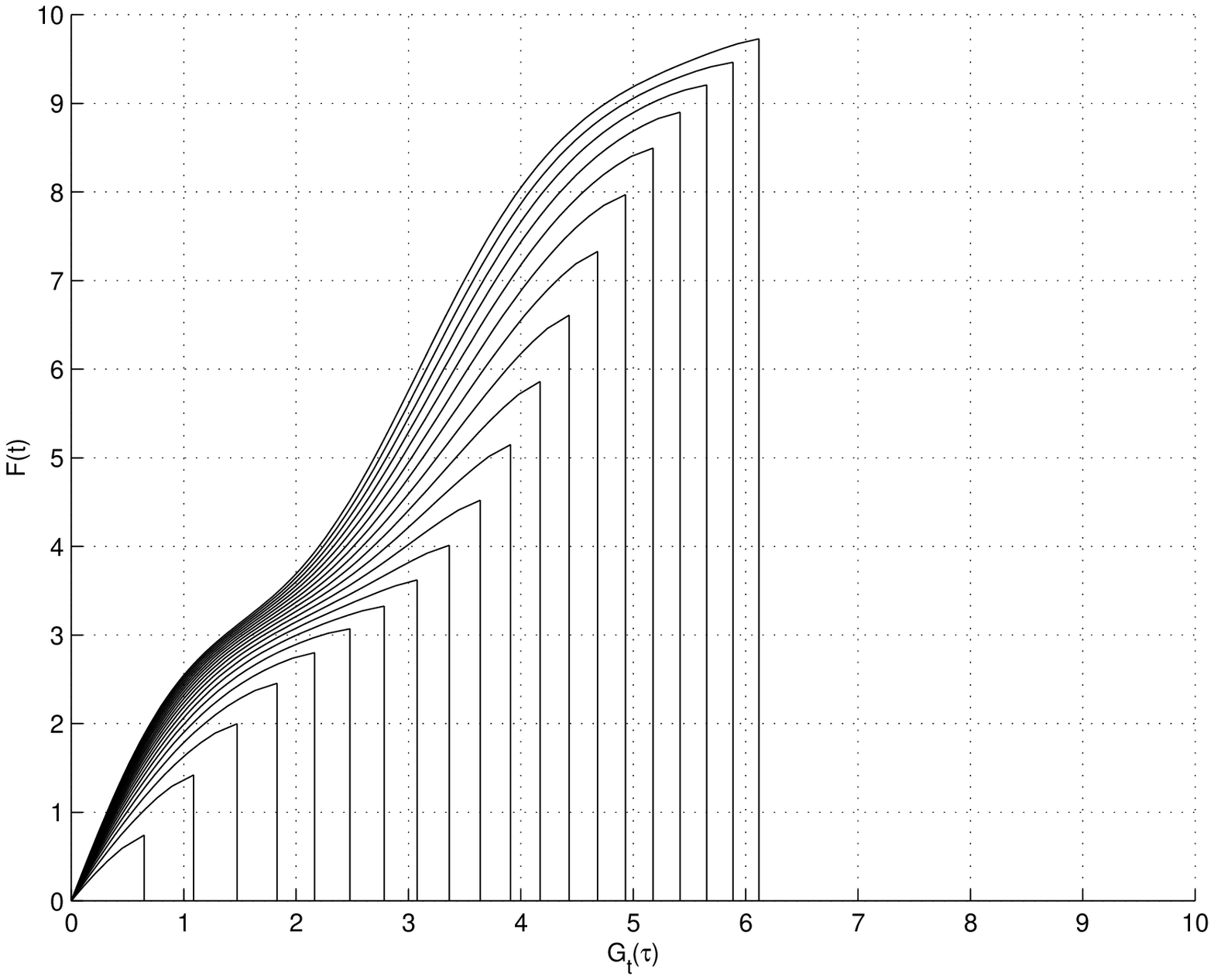}
\caption{\protect\parbox[t]{9cm}{Snapshots of the changing ``shadow'' of changing ``fence''
for  $_{0}I_{t}^{\alpha}f(t)$, $\alpha =0.75$, $f(t)=t+0.5 \sin(t)$,
with the time interval $\Delta t=0.5$ between the snapshops.}}
\label{fig:shadows}
\end{figure}

What happens when $t$ is changing (namely growing)?
As $t$ changes, the ``fence'' changes simultaneously.
Its length and, in a certain sense, its shape changes.
For illustration, see Fig.~\ref{fig:fencebas}.
If we follow the change of the ``shadow'' on the wall $(g,f)$,
which is changing simultaneously with the ``fence''
(see Fig.\ref{fig:shadows}),
then we have a dynamical geometric interpretation of
the fractional integral (\ref{FI}) as a function of $t$.

\subsection[Right-sided Riemann--Liouville fractional integral]
	{Right-sided Riemann--Liouville fractional integral}

Let us consider the right-sided Riemann--Liouville
fractional integral \cite{Podlubny-FDE-book,SKM},
\begin{equation}\label{Right-FI}
	_{t}I_{0}^{\alpha} f(t)
	 =  	\frac{1}{\Gamma (\alpha)}
		\int\limits_{t}^{b}
		f(\tau) (\tau-t)^{\alpha - 1}
		d\tau,
\end{equation}
and write it in the form
\begin{equation}\label{Right-FI-as-SI}
	_{t}I_{0}^{\alpha} f(t) =
		\int\limits_{t}^{b}
		f(\tau)
		dh_{t}(\tau),
\end{equation}
\begin{equation}\label{eq:h-function}
	h_{t}(\tau) = 	\frac{1}{\Gamma (\alpha + 1)}
			\Bigl\{
			t^{\alpha} + (\tau -t)^{\alpha}
			\Bigr\}.
\end{equation}

\begin{figure}
\centering
\includegraphics[width=0.7\textwidth]{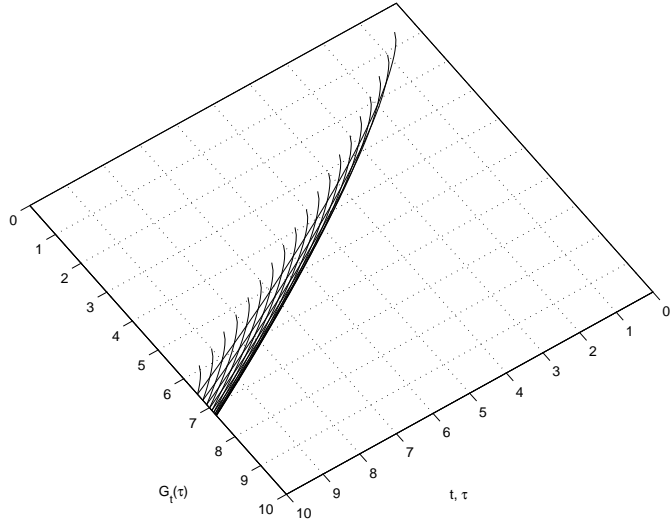}
\caption{\protect\parbox[t]{9cm}{The process of change of the fence basis shape
\protect\\
for $_{t}I_{10}^{\alpha}f(t)$, $\alpha=0.75$.}}
\label{fig:fencebasr}
\end{figure}

Then we can provide a geometric interpretation similar
to the geometric interpretation of the left-sided
Riemann--Liouville fractional integral. However, in this
case there is no any fixed point in the ``fence'' base --
the end, corresponding to $\tau = b$, moves
along the line $\tau=b$ in the plane $(\tau, g)$
when the ``fence'' changes its shape.
This movement can be observed in Fig.~\ref{fig:fencebasr}.
(In the case of the left-sided integral, the left end,
corresponding to $\tau=0$, is fixed and does not move.)

All other parts of the geometric interpretation remain
the same: the ``fence'' changes its shape as $t$ changes
from $0$ to $b$, and the changing shadows of this ``fence''
on the walls $(g,f)$ and $(\tau, f)$ represent correspondingly
the right-sided Riemann--Liouville fractional integral (\ref{Right-FI})
and the classical integral with the moving lower limit:
\begin{equation}
	_{t}I_{b}^{1}(t) =
		\int\limits_{t}^{b}
		f(\tau) d\tau;
\end{equation}

Obviously, for $g_{t}(\tau)=\tau$ both ``shadows'' are equal.
Therefore, we see that not only the left-sided, but also
the right-sided Riemann-Liouville fractional integration
includes the classical definite integration as a particular
case even from the geometrical point of view.

\subsection[Riesz potential]
	{Riesz potential}

The Riesz potential \cite{Podlubny-FDE-book,SKM}
\begin{equation}\label{eq:Riesz-integral}
	_{0}R_{b}^{\alpha} f(t)	=
		\frac{1}{\Gamma (\alpha)}
		\int\limits_{0}^{b}
		f(\tau) |\tau-t|^{\alpha - 1}
		d\tau
\end{equation}
is the sum
of the left-sided and the right-sided
Riemann--Liouville fractional integrals:
\begin{equation}\label{eq:Rp-as-sum-of-LFI-and-RFI}
	_{0}R_{b}^{\alpha} f(t)	=
		\frac{1}{\Gamma (\alpha)}
		\int\limits_{0}^{t}
		f(\tau) (t-\tau)^{\alpha - 1}
		d\tau
		+
		\frac{1}{\Gamma (\alpha)}
		\int\limits_{t}^{b}
		f(\tau) (\tau-t)^{\alpha - 1}
		d\tau.
\end{equation}

The Riesz potential (\ref{eq:Riesz-integral})
can be written in the form
\begin{equation}\label{RP-FI-as-SI}
	_{0}R_{b}^{\alpha} f(t) =
		\int\limits_{0}^{b}
		f(\tau)
		dr_{t}(\tau),
\end{equation}
\begin{equation}\label{eq:r-function}
	r_{t}(\tau) = 	\frac{1}{\Gamma (\alpha + 1)}
			\Bigl\{
			t^{\alpha} + \mbox{sign}{(\tau-t)}\, |\tau -t|^{\alpha}
			\Bigr\}.
\end{equation}

\begin{figure}
\centering
\includegraphics[width=0.7\textwidth]{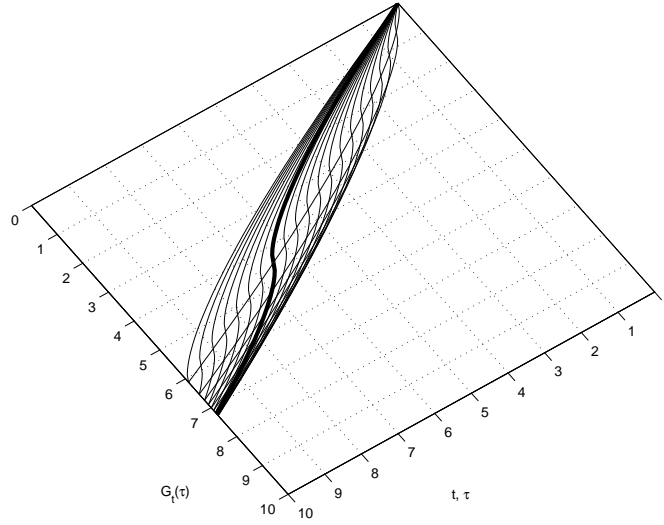}
\caption{\protect\parbox[t]{9cm}{The process of change of the fence basis shape
\protect\\
for the Riesz potential
$_{0}R_{10}^{\alpha}f(t)$, $\alpha=0.75$.}}
\label{fig:rieszbas}
\end{figure}

The shape of the ``fence'', corresponding to the Riesz potential,
is described by the function $r_{t}(\tau)$.
In this case the ``fence'' consists of the two parts:
one of them (for $0<\tau<t$) is the same as in the
case of the left-sided Riemann--Liouville fractional integral,
and the second  (for $t<\tau<b$) is the same as for the right-sided
Riemann--Liouville integral, as shown in Fig.~\ref{fig:rieszbas}.
Both parts are joined smothly at the inflection point $\tau=t$.

The shape of the ``fence'', corresponding to the Riesz potential,
is shown in some of its intermediate position by the bold line in Fig.~\ref{fig:rieszbas}. Obviously, Fig.~\ref{fig:rieszbas}
can be obtained by laying Fig.~\ref{fig:fencebasr} over
Fig.~\ref{fig:fencebas}, which is a geometric interpretation of
the relationship (\ref{eq:Rp-as-sum-of-LFI-and-RFI}).

The shadow of this ``fence''
on the wall $(g,f)$ represents the Riesz potential
(\ref{eq:Riesz-integral}), while the shadow
on the wall $(\tau,f)$ corresponds to the classical
integral
\begin{equation}\label{eq:integral_0_b}
	I(t) = \int\limits_{0}^{b} f(\tau) d\tau.
\end{equation}
For $\alpha=1$ both ``shadows'' are equal.
This shows that the classical definite integral (\ref{eq:integral_0_b})
is a particular case of the Riesz fractional potential
(\ref{eq:Riesz-integral}) even from the geometric point of view.
We have already seen this inclusion in the case of the left-sided
and the right-sided Riemann--Liouville fractional integration.
This demonstrates the strength of the suggested geometric interpretation
of these three types of generalization of the notion of integration.

\subsection{Feller potential}

The Feller potential operator $\Phi^{\alpha}f(t)$
is, similarly to the Riesz potential,
also a linear combination of the left- and right-sided
Riemann--Liouville fractional integrals,
but with general constant coefficients $c,d$
\cite[Chap. 3]{SKM}:
\begin{equation}\label{eq:Feller-potential}
    \Phi^{\alpha} f(t) =
    c \, _{a}I_{t}^{\alpha} f(t) +
    d \, _{t}I_{b}^{\alpha} f(t).
\end{equation}

The geometric interpretation of the Feller potential
can be easily obtained by properly scaling
and then superimposing Fig.~\ref{fig:fencebasr}
and Fig.~\ref{fig:fencebas}. The ``fence'' obtained
in this way is, in general, discontinuous at $\tau=t$.
Its shadow on the wall $(\tau, f)$ is equal to the
classical definite integral (\ref{eq:integral_0_b}).
The shadow on the wall $(g,f)$ consists, in general,
of the two areas, which may overlap depending on the values
of the coefficients $c$ and $d$.

\section{Two kinds of time -- I}

The geometric interpretation of fractional integration,
given in the previous sections, is substantially based
on adding the third dimension (for $g_{t}(\tau)$) to
the classical pair $\tau$, $f(\tau)$.
	If we consider $\tau$ as time, then $g(\tau)$
can be interpreted just as a ``deformed'' time scale.
What could be the meaning of having -- and using -- two
time axes? To answer this question, let us recall
some facts of the history of the development of
the notion of time.

That were contributions of Barrows and Newton to the development
of mathematics and physics in the XVII century which led
to the appearance of the ``mathematical time'', which is
postulated to ``flow equably'' and which is usually
depicted as a semi-infinite straight line \cite{Whitrow}.

Newton himself postulated \cite{Newton}:
\begin{quotation}
\noindent
``Absolute, true and mathematical time of itself, and from its own nature,
flows equably without relation to anything external.''
\end{quotation}

Such a postulate was absolutely necessary for
developing Newton's differential calculus
and applying it to problems of mechanics \cite{Whitrow}:

\begin{quotation}
\noindent
``The outstanding mathematical achievement associated with the
geometrization of time was, of course, the invention of the calculus of fluxions by Newton.''

``Mathematically, Newton seems to have found support for his
belief in absolute time by the need, in principle, for an ideal
rate-measurer.''
\end{quotation}

The invention of differential and integral calculus and today's use
of them is the strongest reason for continuing using homogeneous equably
flowing time.

\begin{figure}
\includegraphics[width=\textwidth]{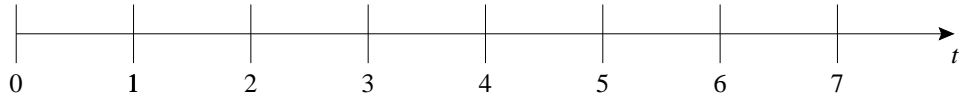}
\caption{Homogeneous time axes.}
\label{fig:homogeneous-time}
\end{figure}

Time is often depicted using the time axis,
and the geometrically equal intervals of the time axis
are considered as corresponding to equal time
intervals (Fig. \ref{fig:homogeneous-time}).

\begin{figure}
\includegraphics[width=\textwidth]{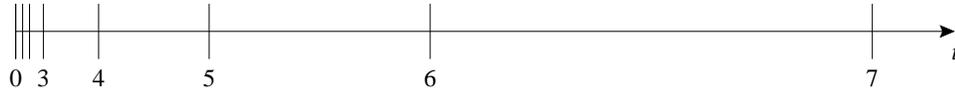}
\caption{Time slowing down.}
\label{fig:2x-time}
\end{figure}

This assumption, however, cannot be neither proved
nor rejected by experiment. Two lengths of geometric
intervals can be measured and compared, since they
are available for measurement simultaneously,
at the same time (or, more precisely, at the same time
and at the same place).
	Two time intervals can never be compared, because
they are available to us for measurement (or for observation)
only sequentially.

Indeed, how do we measure time intervals?
Only by observing some processes, which we
consider as regularly repeated.
G. Clemence wrote \cite{Clemence}:
\begin{quotation}
\noindent
``The measurement of time is essentially a process of counting.
Any recurring phenomenon whatever, the occurences of which
can be counted, is in fact a measure of time.''
\end{quotation}

	Clocks, including atomic clocks,
repeat their ``ticks'', and we simply count
those ticks, calling them hours, minutes,
seconds, milliseconds, etc. But we are not able
to verify if the \textit{absolute} time which elapsed between,
say, the fifth and the sixth tick (the sixth ``second'')
is exactly the same as the time, which  elapsed between
the sixth and the seventh tick (the seventh ``second'').
This possible inhomogeneity of the time scale is illustrated
in Fig.~\ref{fig:2x-time}.

The fact that time measurement as a process of counting of
repeating discrete events does not really exclude
inhomogeneity of time, has been nicely mentioned
by L. Carroll in \emph{Alice's Adventures in Wonderland}
\cite[Chap.~7]{Carroll}:
\begin{quotation}
``\ldots I know I have to beat time when I learn music.''

``Ah! That accounts for it,'' said the Hatter. ``He [Time] won't stand
beating. Now, if you only kept on good terms with him, he'd do almost
anything you liked to do with the clock\ldots ''
\end{quotation}

Figures \ref{fig:homogeneous-time} and \ref{fig:2x-time}
show those ``clock ticks'', which we can register,
only symbolically. One can interpret them as if there
exists some \textit{absolute}, or \textit{cosmic},
inhomogeneous time axis,
to which we can compare \textit{individual} homogeneous time,
represented by some ``clock ticks''.
Our picture of the individual homogeneous time has the form shown
in Fig.~\ref{fig:homogeneous-time}. The cosmic time
may be not necessarily flowing equably, like that
shown in Fig.~\ref{fig:2x-time}.

To illustrate the idea, let us consider the following
situation. Suppose person N has two devices: one is
a speedometer, and another one is the clock, which is slowing
down, so the interval between the two subsequent ticks
is double comparing to the interval between the previous
ticks (see Fig.~\ref{fig:2x-time}). Person N reads
the velocity values indicated by the speedometer at each
encountered ``second'', without knowing that the clock
is, in fact, slowing down.

Using these two series of data, namely the recorded sequence
of values of speed, and the sequence of the counted ``seconds'',
person N can estimate the distance which he has passed.

For simplicity, let us suppose that the first
``second'' of the time shown by the clocks is equal to
the absolute time ``second''.
	The results of observations in this hypothetical
experiment are given in Table 1.

\begin{table}
\centering
\caption{Recording speed using slowing-down clocks}
\begin{tabular}{|c|c|c|}
\hline
\hline
Person N & Recorded values & Observer $O$ \\
individual & of velocity [m/s] & absolute (cosmic) \\
``seconds'' &                  & ``seconds''\\
\hline
\hline
0	&	10	&	0	\\
1	&	11	&	1	\\
2	&	12	&	3	\\
3	&	13	&	7	\\
4	&	12	&	15	\\
5	&	11	&	31	\\
6	&	10	&	63	\\
7	&	 9	&	127	\\
\hline
\hline
\end{tabular}
\end{table}

Person N will compute the distance he has passed as
$$
S_{N}=10 \cdot 1 + 11 \cdot 1 + 12 \cdot 1 + 13 \cdot 1
	+ 12 \cdot 1 + 11 \cdot 1 + 10 \cdot 1 =79.
$$

However, if there would be an independent observer $O$,
knowing about the slowing-down clock,
then such an observer
would obtain a notably different result for the distance
passed by person N:
$$
S_{O}=   10 \cdot 1  + 11 \cdot 2  + 12 \cdot 4 + 13 \cdot 8
       + 12 \cdot 16 + 11 \cdot 32 + 10 \cdot 64 = 1368.
$$

Below we use this idea for giving a new mechanical
interpretation of the Stieltjes integral.


\section{Physical interpretation of \\ the Stieltjes integral}

Imagine a car equipped with two devices for measurements:
the speedometer recording the velocity $v(\tau)$,
and the clock which should show the time $\tau$.
The clock, however, shows the time incorrectly;
let us suppose that the relationship between the wrong time $\tau$,
which is shown by the clock and which the driver considers as
the correct time, on one hand, and the true time $T$, on the other,
is described by the function $T=g(\tau)$.
	This means that where the driver ``measures'' the time
interval $d\tau$, the real time interval is given by $dT=dg(\tau)$.

The driver $A$, who do not know about wrong operation of the clock,
will compute the passed distance as the classical integral:
\begin{equation}
S_{A}(t) = \int\limits_{0}^{t}
	   v(\tau)
	   d\tau \,.
\end{equation}

However, the observer $O$ knowing about the wrong clock and
having the function $g(\tau)$, which restores the correct values
of time from the driver's wrong time $\tau$,
will compute the really passed distance as
\begin{equation} \label{eq:Stieltjes-integral-II}
S_{O}(t)=	\int\limits_{0}^{t}
	v(\tau)
	dg(\tau).
\end{equation}

This example shows that the Stieltjes integral (\ref{eq:Stieltjes-integral-II})
can be interpreted as the real distance passed by a moving object,
for which we have recorded correct values of speed and incorrect
values of time; the relationship between the wrongly recorded time $\tau$
and the correct time $T$ is given by a known function $T=g(\tau)$.


\section{Physical interpretation of fractional integration: Shadows of the past}

Now let us consider the left-sided Riemann--Liouville fractional integral
\begin{equation}\label{eq:S_O}
	S_{O}(t)=
		\int\limits_{0}^{t}
		v(\tau)
		dg_{t}(\tau) = \;\,
_{0}I_{t}^{\alpha} v(t),
\end{equation}
where $g_{t}(\tau)$ is given by (\ref{eq:g-function}).

The fractional integral $S_{O}(t)$ of the function $v(\tau)$
can be interpreted as the real distance passed by a moving object,
for which we have recorded the local values of its speed $v(\tau)$
(individual speed)
and the local values of its time $\tau$ (individual time);
the relationship between the locally recorded
time $\tau$ (which is considered as flowing equably)
and the cosmic time (which flows non-equably)
is given by a known function $g_{t}(\tau)$.

The function $g_{t}(\tau)$ describes the inhomogeneous time scale,
which depends not only on $\tau$, but also on the parameter $t$
representing the last measured value of the individual time
of the moving object.
When $t$ changes, the entire preceding cosmic time interval changes as well.
This is in agreement with the current views in physics.
Indeed, B.~N.~Ivanov \cite[p.~33]{Ivanov} mentioned
that time intervals depend on gravitational fields.
Similarly, S.~Hawking \cite[p.~32--33]{Hawking-BHT}
wrote that:
\begin{quotation}
\noindent
``\dots\ time should appear to run slower near a massive body
like the earth.'' [\dots]

``\dots\  there is no unique absolute time,
but instead each individual has his own personal measure of time
that depends on where he is and how he is moving.''
\end{quotation}

When a moving body changes its position in space--time,
the gravitational field in the entire space--time also changes
due to this movement.
As a consequence, the cosmic time interval, which corresponds
to the history of the movement of the moving object, changes.
This affects the calculation (using formula (\ref{eq:S_O}))
of the real distance $S_{O}(t)$
passed by such a moving object.

In other words, the left-sided Riemann--Liouville fractional integral
of the individual speed $v(\tau)$ of a moving object,
for which the relationship between its individual time $\tau$
and the cosmic time $T$ at each individual time instance $t$
is given by the known function $T=g_{t}(\tau)$ described by
the equation (\ref{eq:g-function}),
represents the real distance $S_{O}(t)$ passed by that object.

\section{Physical interpretation of the Riemann-Liouville \\
         fractional derivative}\label{sec:RL_Phys}

On the other hand, we can
use the properties of fractional differentiation and integration
\cite{Podlubny-FDE-book,SKM}
and express $v(t)$ from the equation (\ref{eq:S_O})
as a left-sided Riemann--Liouville fractional derivative of $S_{O}(t)$:
\begin{equation}\label{eq:v_R-L}
v(t) = \;\, _{0}D_{t}^{\alpha}S_{O}(t)
\end{equation}
where $_{0}D_{t}^{\alpha}$ denotes the Riemann--Liouville fractional
derivative \cite{Podlubny-FDE-book,SKM},
which is for $0<\alpha<1$ defined  by
\begin{equation}
	_{0}D_{t}^{\alpha} f(t) =
	\frac{1}{\Gamma(1-\alpha)}
	\frac{d}{dt}
	\int\limits_{0}^{t}
	\frac{f(\tau)d\tau}
	     {(t-\tau)^{\alpha}}.
\end{equation}

This shows that the left-sided Riemann--Liouville fractional derivative
of the real distance $S_{O}(t)$ passed by a moving object,
for which the relationship between its individual time $\tau$
and the cosmic time $T$ at each individual time instance $t$
is given by the known function $T=g_{t}(\tau)$ described
by equation (\ref{eq:g-function}),
is equal to the individual speed $v(\tau)$ of that object.

On the other hand, we can differentiate
the relationship (\ref{eq:S_O}) with respect
to the cosmic time variable $t$, which gives
the relationship between the velocity $v_{O}(t) = S'_{O}(t)$
of the movement from the viewpoint of the independent
observer $O$ and the individual velocity $v(t)$:
\begin{equation}\label{eq:v_O}
v_{O}(t) = \frac{d}{dt} \; _{0}I_{t}^{\alpha} v(t) =
\; _{0}D_{t}^{1-\alpha} v(t),
\end{equation}

Therefore, the ($1-\alpha$)-th--order Riemann--Liouville derivative
of the individual velocity $v(t)$ is equal to the velocity $v_{O}(t)$
from the viewpoint of the independent observer,
if the individual time $\tau$ and the cosmic time $T$
are related by the function $T=g_{t}(\tau)$ described
by equation (\ref{eq:g-function}).
For $\alpha=1$, when there is no dynamic deformation of the
time scale, both velocities coincide: $v_{O}(t) = v(t)$.

\section{Physical interpretation of the Caputo \\
         fractional derivative}

Applying fractional integration of order $\beta=1-\alpha$ to both parts
of the relationship (\ref{eq:v_O}) gives:
\begin{equation}\label{eq:v_Caputo}
	v(t) = \; _{0}I_{t}^{1-\alpha} v_{O} (t) =
		\; _{0}I_{t}^{1-\alpha} S'_{O} (t) =
		\; _{0}^{C}\! D_{t}^{\alpha} S_{O} (t),
\end{equation}
where $_{0}^{C}\! D_{t}^{\alpha}$ denotes the Caputo fractional
derivative \cite{Caputo-1967,Caputo-ED,Podlubny-FDE-book},
which is for $0<\alpha<1$ defined  by
\begin{equation}
	_{0}^{C}\! D_{t}^{\alpha} f(t) =
	\frac{1}{\Gamma(1-\alpha)}
	\int\limits_{0}^{t}
	\frac{f'(\tau)d\tau}
	     {(t-\tau)^{\alpha}} .
\end{equation}

The relationship (\ref{eq:v_Caputo}) is similar to (\ref{eq:v_R-L}).
Therefore, the Caputo fractional derivative has the same
physical interpretation as the Riemann--Liouville fractional derivative
(see Section~\ref{sec:RL_Phys}).
This coincidence becomes more obvious, if we recall
\cite{Podlubny-FDE-book}
that
if $f(0)=0$, then the Riemann--Liouville derivative and the Caputo derivative
of order $\alpha$ ($0<\alpha<1$),
coincide: $_{0}^{C}D_{t}^{\alpha} f(t) = \;_{0}D_{t}^{\alpha} f(t)$.

\section{Two kinds of time -- II}

The suggested physical interpretation of fractional integration and
fractional differentiation is based on using two kinds of time:
the cosmic time and the individual time.

As mentioned above, due to the history of the development
of mathematics and physics, we are taught to think about the time, in fact,
geometrically. The real roots of this go even far back to Euclid \cite{Whitrow}:
\begin{quotation}
\noindent
``Euclid considered space as the primary concept of science and relegated time
to poor second.''
\end{quotation}

The entire integral and differential calculus is based on using
mathematical (homogeneous, equably flowing) time.
There is no chance to change this state, and there is
nothing to suggest instead of the classical calculus.
Moreover, there is probably even no need for this.
We can just realize that the classical calculus provides tools
for describing the dynamic properties of the cosmic time,
which -- according to physicists -- is inhomogeneous (flowing non-equably).
Indeed \cite[p.~33--34]{Hawking-BHT},
\begin{quotation}
\noindent
``The old idea of an essentually unchanging universe that could
have existed, and could continue to exist, forever was replaced
by the notion of a dynamic, expanding universe that seemed to have
begun a finite time ago\dots''
\end{quotation}

Clearly, the expansion of the universe implies that neither spatial
scale nor time scale remains homogeneous; they both are dynamic.
	For describing the inhomogeneous time, the ideal homogeneous
time scale can be used.
This approach is not new; it has already been
used in the theory of relativity for describing shortening of time
intervals. This means that in fact two time scales are considered
simultaneously: the ideal, equably flowing homogeneous time, and the
cosmic (inhomogeneous) time. The change of scale of the cosmic time is
described using the homogeneous time scale as a reference scale.
In other words, the homogeneous time scale is just an ideal notion,
which is necessary for developing mathematical models
describing inhomogeneous cosmic time and its change.
In this respect we can, without discussing other views on this subject,
recall the remark made by A. Daigneault and A. Sangalli
in their essay \cite{Daigneault-Sangalli}
about I.~E.~Segal and his two-time cosmology
(``chronometric cosmology'', or CC) \cite{Segal} --
note that ``\textit{perhaps!}'':
\begin{quotation}
\noindent
``According to CC, Einstein's model is the correct one to understand
the universe as a whole (i.e., global space--time), except that there
are two kinds of time: a cosmic or Einstein's time $t$, and a local
or Minkowski's time $x_0$, which is (perhaps!) the time measured
by existing techniques. [\ldots] Simply put, Einstein's cosmic time
is the ``real'' one, whereas Minkowski's time is only an approximation
of~$t$.''
\end{quotation}

\noindent
So, the ideal model of equably flowing homogeneous time
can be considered as a rough approximation of the cosmic time.

\section{Geometric and physical interpretation \\
         of the Volterra convolution integral}

It should be mentioned that we can also provide a geometric
and physical interpretation for more general integrals.

The Riemann--Liouville fractional integral is a particular
case of convolution integrals of the Volterra type:
\begin{equation}
	K * f(t) =
	\int\limits_{0}^{t}
	f(\tau)
	k(t-\tau)
	d\tau
\end{equation}

Assuming that $k(t)=K'(t)$, we can write this integral in the form
\begin{equation}
	K * f(t) =
	\int\limits_{0}^{t}
	f(\tau)
	dq_{t}(\tau),
\end{equation}
\begin{equation}
	q_{t}(\tau) = K(t) - K(t-\tau).
\end{equation}

The geometric and physical interpretation of the Volterra convolution
integral is then similar to the suggested interpretations
for fractional integrals. The function $q_{t}(\tau)$
determines the changing shape of the ``live fence''
(in the case of the geometric interpretation,
see Figs.~\ref{fig:fence} and \ref{fig:fencebas})
and the relationship between the individual time
and the cosmic time of a moving object
(in the case of the physical interpretation).

\section*{Acknowledgment}

This work has been done during the author's visit
to Department of Applied Mathematics and Theoretical Physics (DAMTP)
of the University of Cambridge, UK, in July 2000, supported
by the Cambridge Colleges Hospitality Scheme.

The author will always
remember late Professor David G. Crighton, the former Head of DAMTP
and the Master of the Jesus College,
who invited the author and arranged all details
of the author's stay and research at DAMTP.
Without this the present paper could not appear.

The author also expresses his deep gratitude to Professor
Stephen Heath, the President of the Jesus College,
and to all members of the College, for their cordial approach
during the whole period of the author's visit.


\end{document}